\documentclass[12pt]{article}
\usepackage{graphicx}
\usepackage{amsmath}
\usepackage{amssymb}

\title{Periodic auxetics: Structure and design}
\author{Ciprian S. Borcea\footnote{Department of Mathematics, Rider University,  Lawrenceville, NJ, USA.
 \texttt{borcea@rider.edu} }  \ and Ileana Streinu\footnote{Computer Science Department, Smith College, Northampton, MA, USA. \texttt{istreinu@smith.edu}}}

\date{}

\begin{document}
\maketitle

\begin{abstract}
Materials science has adopted the term of auxetic behavior for structural deformations where stretching in some direction entails lateral widening, rather than lateral shrinking. Most studies, in the last three decades, have explored repetitive or cellular structures and used the notion of negative Poisson's ratio as the hallmark of auxetic behavior. However, no general auxetic principle has been established from this perspective. In the present  paper,  we show that a purely geometric approach to periodic auxetics is apt to identify essential characteristics of frameworks with auxetic deformations and can generate a systematic and endless series of periodic auxetic designs. The critical features refer to convexity properties expressed through families of homothetic ellipsoids.
\end{abstract}


\section*{Introduction}

A flexible structure exhibits auxetic behavior when stretching in a given direction involves widening in orthogonal directions. In elasticity theory, this type of behavior amounts to {\em negative Poisson's ratios} (for orthogonal directions paired with the given tensile direction). Materials and structures with negative Poisson's ratios have attracted increased attention after the 1987 publication of Lakes' results on foams \cite{lakes:negativePoisson:1987}. The term {\em auxetic} (from {\em increase, growth}) gained currency after \cite{evans:etAl:molecularNetwork:Nature:1991}. The extensive literature generated in  three decades of auxetic studies is reviewed in several recent surveys e.g. \cite{greaves:lakes:etAl:PoissonRatio:2011,lee:singer:thomas:microNanoMaterials:advMat:2012,elipe:lantada:auxeticGeometries:2012,huang:chen:negativePoisson:2016,saxena:auxeticsResearch:advancedEngr:2016,park:auxetic:appliedPhysics:2016,kolken:auxetic:rscAdvances:2017,bertoldi:cellularReview:2017,lakes:auxeticReview:2017}. 

\medskip
Although the role of the underlying geometry is widely recognized, the approach via Poisson's ratios involves modeling assumptions about and dependence on physical properties of the material under consideration. The catalog of auxetic structures found in this literature is rather confined and no general principles of auxetic design have been formulated.

\medskip
The main purpose of the present paper is to show that a {\em purely geometric} approach to periodic auxetics, introduced by the authors in \cite{borcea:streinu:GeomAuxetics:RSPA:2015}, leads to new insights and solves two fundamental problems on structure and design, namely, what geometry is required for
auxetic behavior and how to construct periodic frameworks with auxetic deformations.

\section{Distinctive constituents of the geometric theory}
\label{sec:intro}

Our theory applies to {\em periodic bar-and-joint frameworks}. The dimensions most important for applications are two and three, but the theory works in arbitrary dimension $d$. Periodic frameworks represent the next of kin to finite linkages.

\medskip
The foundations of a deformation theory for periodic frameworks were presented in
	\cite{borcea:streinu:PeriodicFF:RSPA:2010}, with rigidity aspects developed in \cite{borcea:streinu:MinimallyRigid:BLMS:2011} and additional perspectives discussed in \cite{borcea:streinu:CrystallographicSym:PhT:2014,borcea:streinu:Diamond:imaMath:2017}. In \cite{borcea:streinu:LiftStress:DCG:2015}, we formulated and proved a periodic analog of Maxwell's theorem on liftings and stresses and used it to obtain a complete characterization of {\em expansive behavior} in dimension two through the notion of periodic pseudotriangulation \cite{borcea:streinu:kinematicsExpansive:ark14:2014,borcea:streinu:LiftStress:DCG:2015}. A pseudotriangle is a simple planar polygon with exactly three vertices on its convex hull. {\em Expansive behavior}, when all distances between
pairs of vertices increase or stay the same, clearly satisfies the `growth' requirements of {\em auxetic behavior}, that is, {\em expansive implies auxetic}.

\medskip \noindent
{\bf Periodic auxetics.}\ The theory introduced in \cite{borcea:streinu:GeomAuxetics:RSPA:2015} looks {\em directly} at deformation trajectories which exhibit auxeticity. The `lateral widening upon stretching', or rather, in the reverse sense, `lateral shrinking upon compression' is adequately captured in the notion of {\em contraction} (norm at most one) for the linear transformation which takes the periodicity lattice at one moment to the periodicity lattice at a subsequent moment. 

\medskip \noindent
This characterization is completely geometric, valid in arbitrary dimension $d$ and focused on the evolution of the periodicity lattice of the framework. After a choice of independent generators, the deformation of the periodicity lattice gives  a curve in the space of symmetric $d\times d$ matrices, through the variation of the Gram matrix of these generators. In this language, a one-parameter deformation is auxetic when all velocity vectors along the curve belong to the {\em positive semidefinite cone} and {\em strictly} auxetic, when velocity vectors are in the {\em positive definite cone}.

\medskip \noindent
{\bf Remarks.}\
In contrast to the conventional route via Poisson's ratios, this
geometric approach is formulated directly and exclusively in terms of structure and periodicity. The mathematical model of periodic frameworks is in agreement with the rigid unit mode theory used for crystalline materials \cite{dove:displacive:1997}. The bar-and-joint mechanical model is frequently adopted in the literature \cite{evans:auxetic:actaMetal:1994,grima:negativePoisson:RSPA:2012,mitschkeEtAl:geometryLeading:2016} and is gaining critical importance in {\em structural design} due to recent advances in digital manufacturing and the quest for mechanical metamaterials \cite{reisEtAl:designerMatter:2015}.

\medskip 
Besides being completely rigorous, the geometric approach brings the following
advantages:

\medskip
(i) a linear {\em infinitesimal theory} which allows detection of auxetic
capabilities via semidefinite programming \cite{borcea:streinu:GeomAuxetics:RSPA:2015} or other algorithms \cite{borcea:streinu:AuxeticElliptic:2016};

\medskip
(ii) for several degrees of freedom, the notion and description of the {\em infinitesimal auxetic cone}, which is a {\em spectrahedral cone} 
(i.e. a section of the positive semidefinite cone by a linear subspace)  or a linear preimage of a spectrahedral cone
\cite{borcea:streinu:GeomAuxetics:RSPA:2015};

\medskip
(iii) a theorem of isomorphism for auxetic cones under affine transformations of periodic frameworks \cite{borcea:streinu:NewPrinciples:SIAGA:2017};

\medskip
(iv) a general method for converting finite linkages with adequate deformations into periodic frameworks with auxetic deformations \cite{borcea:streinu:NewPrinciples:SIAGA:2017}. This method provides {\em infinite series} of auxetic designs.

\medskip 
As mentioned earlier, the literature based on the traditional approach has generated only a limited list of sporadic patterns and considers the invention of new auxetic designs a {\em challenging problem} \cite{lee:singer:thomas:microNanoMaterials:advMat:2012}, p. 4792.

\section{Outline of the main results}

In this paper, we give comprehensive answers to the following central problems:

\medskip
(1)\ What structural characteristics are required for auxetic capabilities?

%\medskip
(2)\ How to design periodic frameworks with auxetic deformations?

\medskip 
We give first a description of key facts and ideas. Technical details are provided in the next section.

\medskip \noindent
{\bf (1)\ Structural characteristics.}\ The rather stringent structural constraints
which must be satisfied by a periodic framework with auxetic deformations become
visible when using {\em lattice coordinates} i.e. coordinates relative to the periodicity lattice, with a chosen basis of generators \cite{borcea:streinu:CrystallographicSym:PhT:2014}. As detailed below, these coordinates implicate directly the Gram matrix of the generators and the existence of a strictly auxetic infinitesimal deformation has a {\em geometric formulation} which brings forth a {\em family of homothetic ellipsoids}. This family is indexed
by all pairs of vertex orbits with at least one edge orbit connection. All edge vector representatives emanating from a vertex $v_i$ to the orbit of, say, $v_j$ must have their origin and endpoint on the associated ellipsoid. In particular, these points must be in {\em strictly convex position}, a
fairly restrictive necessary condition (which, in dimension three, is visually
easy to assess).
There are also {\em cycle vanishing relations} involving the centers of these ellipsoids. 

\medskip \noindent
The essential fact is that a strictly auxetic infinitesimal deformation is 
completely expressed through a {\em geometric diagram with homothetic ellipsoids}.

\begin{figure}[h]
 \centering
 {\includegraphics[width=0.75\textwidth]{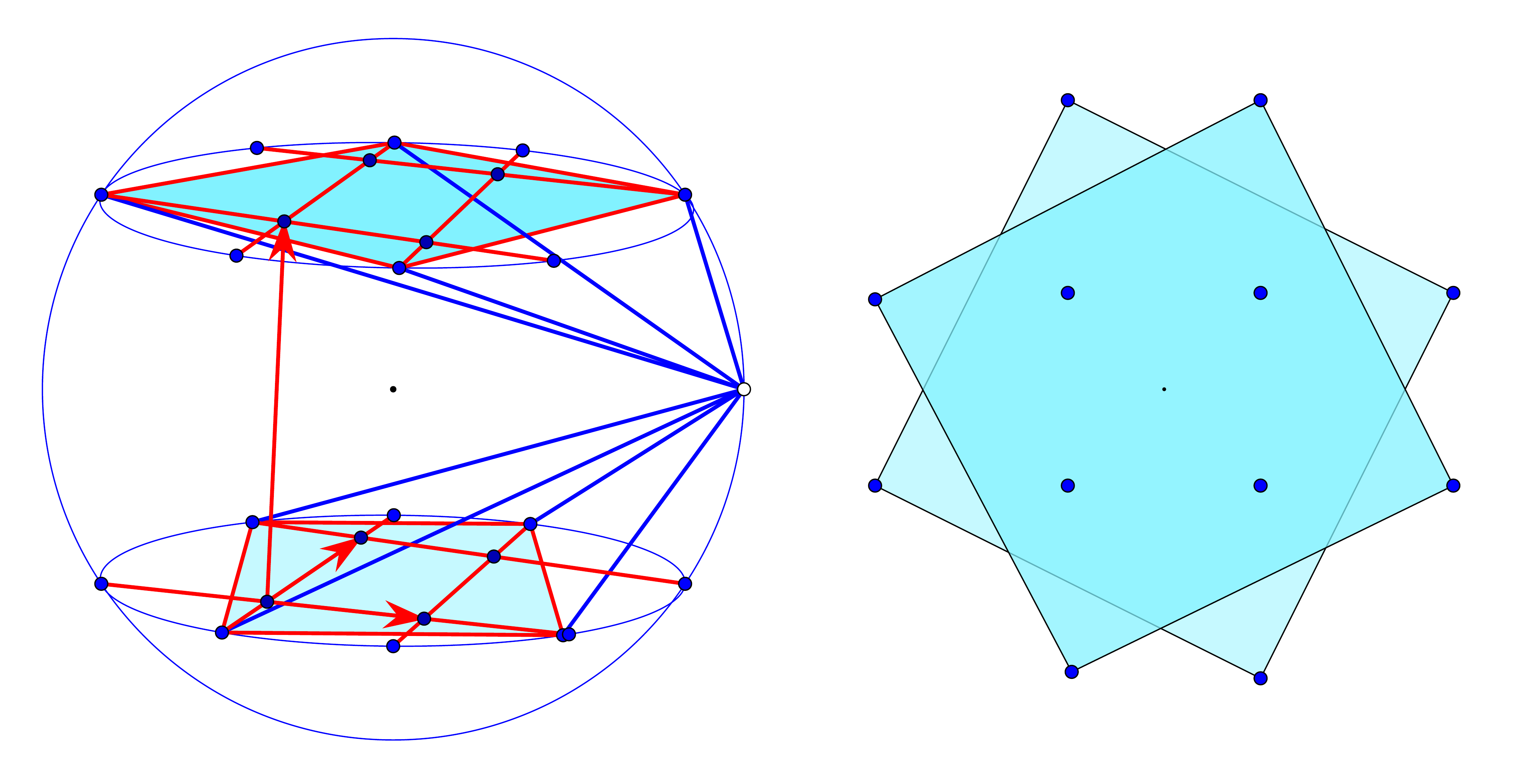}}
 \caption{A blueprint for a three-dimensional auxetic periodic framework with one degree of freedom.
Periodicity generators are shown as red arrows. The framework has two vertex orbits and eight edge orbits. The white vertex and the endpoints of the eight bars emanating from it are on a sphere.
The vertical projection explains why all vectors between vertices in solid blue are in the lattice spanned by the generators.}
 \label{Cover}
\end{figure}

\medskip \noindent
{\bf (2)\ Auxetic design.}\ The design methodology is based on the possibility to satisfy the
requirements identified above in a systematic way. One condition needs particular attention, namely, that vectors which must be periods should generate (over $Z$) a lattice. In Figure~\ref{Cover}, we illustrate
a case with two vertex orbits, respecting this condition. In general, the
lattice condition will be satisfied when we construct our designs using 
{\em rational points} as detailed below.

The design procedure starts with a given quotient graph and can generate an {\em infinite series} of periodic frameworks (with that quotient). By construction,
the generated frameworks will have a strictly auxetic infinitesimal deformation.
A guaranteed one-parameter deformation (extending the infinitesimal one) is
obtained whenever the local deformation space is smooth. In this respect, our
results in \cite{borcea:streinu:MinimallyRigid:BLMS:2011} can be used for selecting quotient graphs with
generic liftings which do have smooth local deformation spaces. The number of
degrees of freedom can be predetermined as well.

Thus, our design method can generate an \textit{infinite virtual catalog}
of periodic frameworks with auxetic deformations. While our paper \cite{borcea:streinu:NewPrinciples:SIAGA:2017}
has already demonstrated that infinite series of auxetic designs are possible
in any dimension $d\geq 2$, the new method is \textit{comprehensive} and amenable
to algorithmic treatment.

\section{Technical details, statements and proofs}

\medskip \noindent
{\bf Periodic graphs and (reduced) quotient graphs.} 
A $d$-periodic graph is a pair $(G,\Gamma)$, where $G=(V,E)$ is a simple infinite graph with vertices $V$, edges
$E$ and finite degree at every vertex, and $\Gamma \subset Aut(G)$ is a free Abelian group of automorphisms
which has rank $d$, acts without fixed points and has a finite number of vertex (and hence, also edge) orbits.
We assume $G$ to be connected. The group  $\Gamma$ is thus isomorphic to $Z^d$ and is called the {\em periodicity group}  of the periodic graph $G$. Its elements $\gamma \in \Gamma \simeq Z^d$ are referred to as {\em periods} of $G$. 
We denote by $G/\Gamma = (V/\Gamma, E/\Gamma)$ the {\em quotient (multi)graph}, with $n=|V/\Gamma|$ the number of vertex orbits and $m=|E/\Gamma|$ the number of edge orbits. We also make use of an associated {\em reduced} quotient graph, where multiple edges have been replaced with simple edges.

\medskip \noindent
{\bf Periodic frameworks and deformations.} 
A periodic placement (or simply placement) of a $d$-periodic graph $(G,\Gamma)$ in $ \mathbb{R}^d$ is defined by two functions:
$p:V\rightarrow  \mathbb{R}^d$ and  $\pi: \Gamma \hookrightarrow {\cal T}( \mathbb{R}^d)$, 
where $p$ assigns points in $ \mathbb{R}^d$ to the vertices $V$ of $G$ and $\pi$ is a faithful representation of the periodicity group $\Gamma$, that is, an injective homomorphism of $\Gamma$ into the group ${\cal T}( \mathbb{R}^d)$ of translations in the Euclidean space $ \mathbb{R}^d$, with $\pi(\Gamma)$ being a lattice of rank $d$. These two functions must satisfy the natural compatibility condition $p(\gamma v)=\pi(\gamma)(p(v))$.

\medskip \noindent
A periodic framework ${\cal F}=(G,\Gamma,p,\pi)$ is a periodic graph $(G,\Gamma)$ together with a placement $(p, \pi)$. Edges are represented as straight segments 
between their endpoints. Two frameworks are considered equivalent when one is obtained from the other by a Euclidean isometry. A {\em one-parameter deformation of the periodic framework} ${\cal F}$ is a  (smooth) family of placements  $p_{\tau}: V\rightarrow  \mathbb{R}^d$  parametrized by time $\tau \in (-\epsilon, \epsilon)$ in a small neighborhood of the initial placement $p_0=p$, which satisfies two conditions: (a) it maintains the lengths of all the edges $e\in E$, and (b) it maintains periodicity under $\Gamma$, via faithful representations $\pi_{\tau}:\Gamma \rightarrow {\textbf T}( \mathbb{R}^d)$ which {\em may change with $\tau$ and give an associated variation of the periodicity lattice of translations} $\pi_{\tau}(\Gamma)$. 
 
\medskip \noindent
After choosing an independent set of $d$ generators for the periodicity lattice $\Gamma$, the image $\pi(\Gamma)$ is completely described via the $d\times d$ matrix $\Lambda$ with column vectors $(\lambda_i)_{i=1,\cdots,d}$ given by the images of the generators under $\pi$. 
The Gram matrix for this basis will be $\omega=\Lambda^t\cdot \Lambda$.

\medskip \noindent
{\bf Infinitesimal deformations in lattice coordinates.} 
Let us fix now a complete set of vertex representatives $v_0,v_1,...,v_{n-1}$
for the $n$ vertex orbits of $(G,\Gamma)$. The framework ${\cal F}$ has them
positioned at $p_i=p(v_i)$. When we pass from these Cartesian coordinates to
lattice coordinates $q_i$, we consider $v_0$ to be the origin, that is $q_0=0$,
and then $\Lambda q_i=p_i-p_0$.

\medskip \noindent
We recall the form of the equations expressing the constant (squared) length of edges, when using parameters $(q_1,...,q_{n-1},\omega)$, cf. \cite{borcea:streinu:CrystallographicSym:PhT:2014}, formula (4.1). Let us consider an edge (denoted here simply by $e_{ij}$) which goes from $v_i$ to a vertex in the orbit of $v_j$. Then, in {\em Cartesian coordinates}, the edge vector is given by $p_j+\lambda_{ij}-p_i$, with some period $\lambda_{ij}=\Lambda n_{ij}\in \pi(\Gamma)$ and $n_{ij}\in \mathbb{Z}^d$. In {\em lattice coordinates}, the edge vector is given by $e_{ij}=q_j+n_{ij}-q_i$ and the squared-length equation is:

\begin{equation}\label{sqLength}
	\ell(e_{ij})^2	=\langle \omega e_{ij},e_{ij} \rangle
\end{equation}

\medskip \noindent
For infinitesimal deformations we use the notation $(\dot{q}_1,...,\dot{q}_{n-1},\dot{\omega})$. From (\ref{sqLength}) we obtain:

\begin{equation}\label{infDef}
 \langle \dot{\omega} e_{ij}, e_{ij} \rangle + 2\langle \omega e_{ij}, \dot{e}_{ij}\rangle = 0 
\end{equation}

\noindent
and since $\dot{e}_{ij}=\dot{q}_j-\dot{q}_i$, we have
 
\begin{equation}\label{quadric} 
	\langle \dot{\omega} e_{ij}, e_{ij} \rangle + 
2\langle \omega (\dot{q}_j-\dot{q}_i), e_{ij} \rangle = 0 
\end{equation}

\medskip \noindent
{\bf Key observation.}\ The definition of auxetic behavior from \cite{borcea:streinu:GeomAuxetics:RSPA:2015},  reviewed in Section \ref{sec:intro}, implies that an infinitesimal deformation is {\em strictly auxetic} precisely when $\dot{\omega}$ is {\em positive definite}. In this case, the geometric meaning of equation (\ref{quadric}) is that {\em all edge vectors} $e_{ij}$ from the orbit of $v_i$ to the orbit of $v_j$, when seen as fixed vectors, have all their endpoints on the
{\bf ellipsoid}:
\begin{equation}\label{ellipsoid} 
	\langle \dot{\omega} x, x \rangle + 
2\langle \omega (\dot{q}_j-\dot{q}_i), x \rangle = 0 
\end{equation}

\noindent
which passes through their common origin $x=0$ and has its center $c_{ij}$
at:
\begin{equation}\label{center}
	c_{ij}= - \dot{\omega}^{-1}\omega (\dot{q}_j-\dot{q}_i)=
	       \dot{\omega}^{-1}\omega (\dot{q}_i-\dot{q}_j)
\end{equation}

\medskip \noindent
{\bf Homothetic ellipsoids on the reduced quotient graph.} Since all the ellipsoid formulae in (\ref{ellipsoid}) have the same quadratic part $\langle \dot{\omega} x, x \rangle$,  they are all {\em homothetic} and can be
labeled by the edges $ij$ of the reduced quotient graph of the given framework.   
Under the assumption of a strictly auxetic infinitesimal deformation, it follows immediately that the quotient graph has no loops. 
	
\medskip \noindent
{\bf Conditions on cycles.}\ Whenever we have a {\em cycle} in the reduced
quotient graph, the sum of the vectors $c_{ij}$ is zero (with the ordering $ij$ corresponding to traversing the cycle). In other words, we obtain (up to translation) a {\em placement of the reduced quotient graph}, with edge vectors
$c_{ij}$ (or $2c_{ij}$).

\medskip \noindent
It will be observed that, while these properties have been discovered in lattice coordinates, they hold true in Cartesian coordinates as well, since the two representations are related by an affine transformation, which takes homothetic
ellipsoids to homothetic ellipsoids and zero-sum vector relations to zero-sum vector relations.

\begin{figure}[h]
 \centering
 {\includegraphics[width=0.65\textwidth]{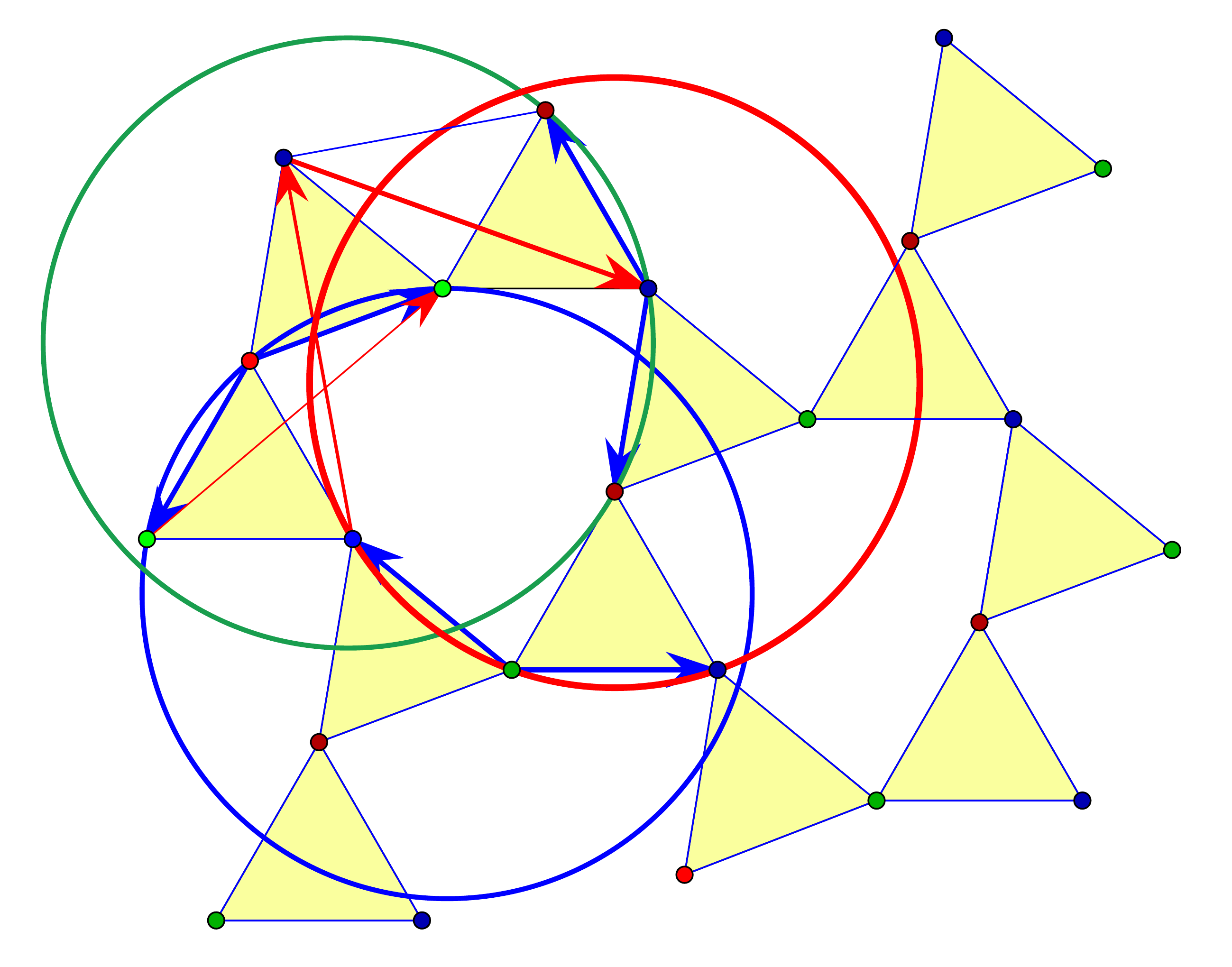}}
 \caption{A Kagome framework has one degree of freedom and its deformation is strictly auxetic until the unit cell reaches maximal area. The geometric manifestation of this fact is the existence of three homothetic ellipses (which are circles, as a result of symmetry) satisfying the cycle condition.}
 \label{FigKagome}
\end{figure}

\medskip \noindent
{\bf Illustration with Kagome frameworks.}\ The structural features disclosed above
will be necessary present in planar frameworks of Kagome type. Such frameworks have
$n=3$ vertex orbits, $m=6$ edge orbits and their deformation space is one dimensional. The quotient graph is a triangle with doubled edges. The auxetic character of Kagome periodic mechanisms with congruent equilateral triangles has been noticed earlier \cite{lubensky:twistedKagome:2012,borcea:streinu:LiftStress:DCG:2015}. 
In Figure~\ref{FigKagome}, the three vertex orbits are colored in red, blue and green. The three ellipses are
actually circles. The red circle passes through a green vertex and the two blue vertices connected with it. The centers of the circles are not shown, but the cycle condition is clearly satisfied by symmetry.

\medskip
\noindent
{\bf Theorem.}\
The reduced quotient graph (with its placement) and the family of homothetic ellipsoids (\ref{ellipsoid}),
with all the edge vector representatives in them allows the reconstruction
of the periodic framework and (up to a scalar) the corresponding strictly auxetic 
infinitesimal deformation.

\medskip \noindent
{\em Proof:}\ The periodic framework is determined by the complete set of 
edge vector representatives. The strictly auxetic infinitesimal line
$(\dot{q}_1:...:\dot{q}_{n-1}:\dot{\omega})$ is retrieved since we have $\dot{\omega}$, up to a scalar, from any of the ellipsoids, and then the 
linear system (\ref{center}) for $\dot{q}_i$ is compatible due to the 
vanishing of sums over cycles observed above.

\medskip \noindent
{\bf Simplifying observation.}\ Since affine transformations of a periodic framework preserve infinitesimal characteristics (such as strict auxeticity),
we may adopt coordinates which turn all our homothetic ellipsoids into spheres.

\medskip \noindent
{\bf Generating periodic frameworks with auxetic capabilities.}\ 
We can construct (generate) such frameworks by following {\em in reverse} the steps of the previously described procedure. Particular attention must be paid to the vectors between endpoints of edge vector representatives in any given ellipsoid and to sums of edge vectors over cycles in the reduced quotient graph: these vectors must be periods and must generate, by linear combinations with  {\em integer coefficients}, a (periodicity) lattice.

\medskip \noindent
This {\em rational dependence} requirement is always fulfilled if we operate only with rational points (i.e. points with rational coordinates). More precisely, we have the following scenario  (in arbitrary dimension $d\geq 2$).

\medskip \noindent
{\bf Construction procedure.}\ We start with 
a (connected) finite multi-graph without loops 
(which plays the role of a given quotient graph). 
We want to construct periodic frameworks with 
strictly auxetic (infinitesimal) capabilities 
and with the given quotient graph.

\medskip  
(1)	We take the reduced graph of the given multi-graph.

% \medskip
(2) We choose a (general, non-degenerate) placement of its vertices in $Q^d$ i.e. at rational points in $ \mathbb{R}^d$ and represent the edges by segments between adjacent vertices.

% \medskip
(3) We draw all spheres with these edges as diameters (hence centers are also rational points). 
	
\medskip \noindent
{\bf Lemma.}
Rational points are dense on all these spheres.

\noindent
{\em Proof:}\ The standard stereographic projection (from a point $e_1$ on the unit sphere $S^{d-1}$  to the tangent space $T_{-e_1}(S^{d-1})$ at the opposite point $-e_1$) gives a bijection between the
rational points on the unit sphere and rational points in the tangent space $T_{-e_1}(S^{d-1})$. Thus, rational points are dense in the unit sphere. It follows that rational points
are dense on any sphere with rational center and at least one rational point.

 \medskip 
(4) When $ij$ is an edge, we choose as many rational points on the corresponding sphere as dictated by the multiplicity in the given multi-graph. (Again, the choice in supposed to be non-degenerate.) The vectors from
vertex $v_i$ to the chosen points will represent all edge vectors connecting
the orbit of $v_i$ to the orbit of $v_j$. 

(5) If degenerate choices are avoided, the span of all vectors between chosen points in any given sphere and all sums of
edge vectors over cycles in the reduced quotient graph will be a (rank $d$) lattice, the periodicity lattice of the framework generated by articulating
the chosen edge representatives.

\begin{figure}[h]
\centering
 {\includegraphics[width=0.85\textwidth]{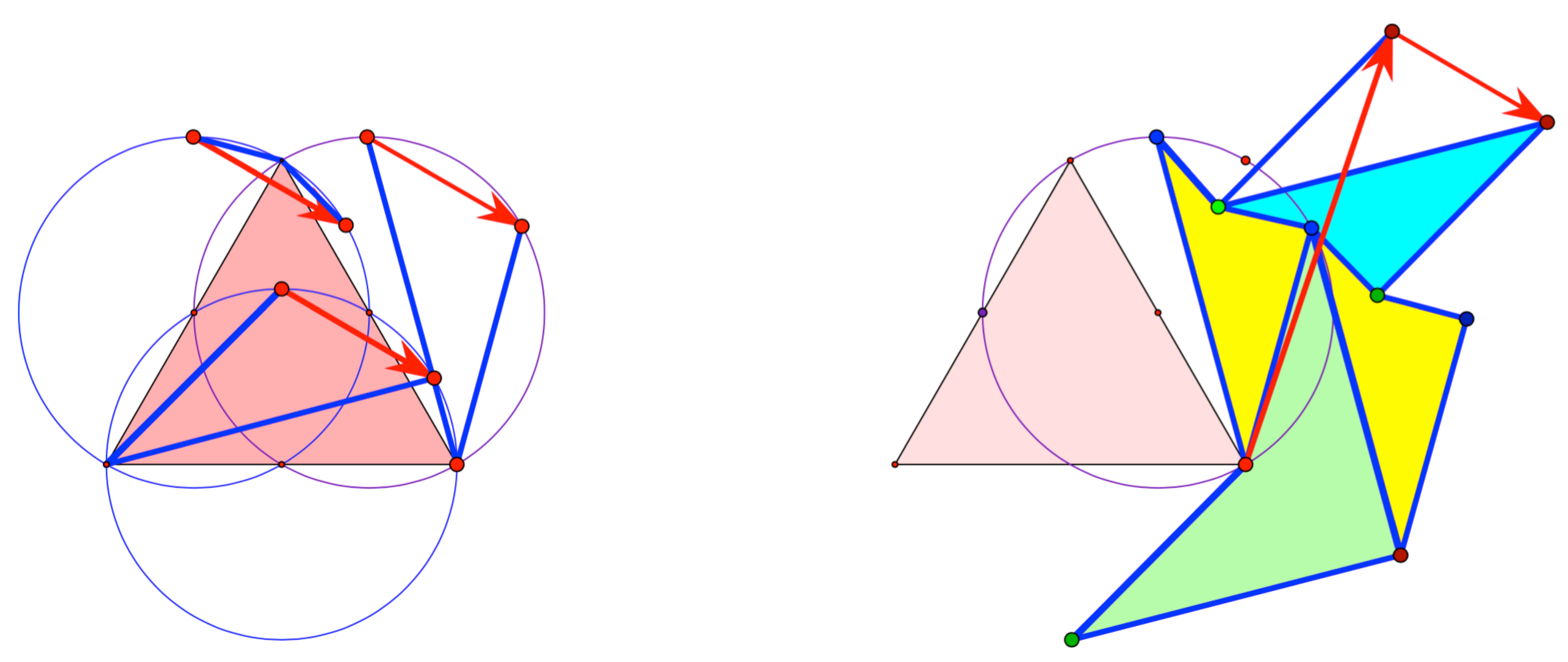}}
 \caption{A planar example starting with a multigraph with $n=3$ vertices and double edges between any pair of vertices ($m=6$). The reduced graph is represented by the highlighted equilateral triangle. The period in each circle is the same. The framework generating process is shown on the right. In this case, the result is a periodic pseudotriangulation.}
 \label{FigPseudo}
\end{figure}

\medskip \noindent
{\bf An explicit example.}\ To illustrate the construction, we use the planar case shown in Figure~\ref{FigPseudo}.
The quotient multi-graph for the example (not shown) has three vertices and double edges between any pair
of vertices. The reduced graph is placed as an equilateral triangle in the plane, and three circles are drawn with diameters given by the three edges. On each circle, we have to position two edge vectors, for the two edge orbits in the multi-graph corresponding to its diameter edge. Choosing a counter-clockwise orientation for the triangle edges, the edge vectors emanate from the source vertex of the corresponding reduced oriented edge. The left hand side of the figure shows this stage of the construction, with the added particularity that in each circle, the red vector between the endpoints of the two edges is the same free vector. This particular choice avoids the need to maneuver with rational points, since this periodicity vector, together with one resulting from the single cycle in the
reduced graph, will span the periodicity lattice. 

The periodic framework on the right is now generated from the diagram constructed on the left, as follows. We start at one vertex and continue to place (all allowed) edge vectors emanating from the endpoints of the initial two edges. Indefinite repetitions of this type of edge addition will generate the entire (connected) periodic framework. The right hand side of the figure shows the stage reached after adding enough edges to see (in green, yellow and blue) representatives of the three type of pseudotriangular faces of the resulting planar framework. Obviously, the initial triangle and the circle are not part of the framework. As a periodic pseudotriangulation, the framework has one degree of freedom and its local deformation is not only auxetic, but actually expansive \cite{borcea:streinu:LiftStress:DCG:2015}.

\medskip
\noindent
{\bf An example with higher symmetry.}\ 
We discuss now a three-dimensional example based on a multi-graph with four vertices and double edges between any pair of vertices. The reduced graph is
represented as a regular tetrahedron with vertices at: $(-1,-1,-1)$, $(1,1,-1)$,
$(1,-1,1)$ and $(-1,1,1)$. We trace the six spheres with diameters given by
the six edges of this tetrahedron. Our choices of edge vector representatives for
the framework will use the symmetry group of the regular tetrahedron, represented
in our setting by permutations of the three coordinates and an even number of sign changes in front of them. Thus, it will be enough to indicate 
one of the two edge vectors
emanating from $(-1,-1,-1)$ and ending at vertices in the orbit represented by
$(1,1,-1)$. We choose it to be of the form $(\alpha,\beta,0)$. The condition that,
when emanating from $(-1,-1,-1)$, the endpoint of the vector must be on the sphere
over the tetrahedral edge from $(-1,-1,-1)$ to $(1,1,-1)$ amounts to: $(\alpha-1)^2+(\beta-1)^2 = 2$.

\begin{figure}[t]
\centering
 {\includegraphics[width=0.40\textwidth]{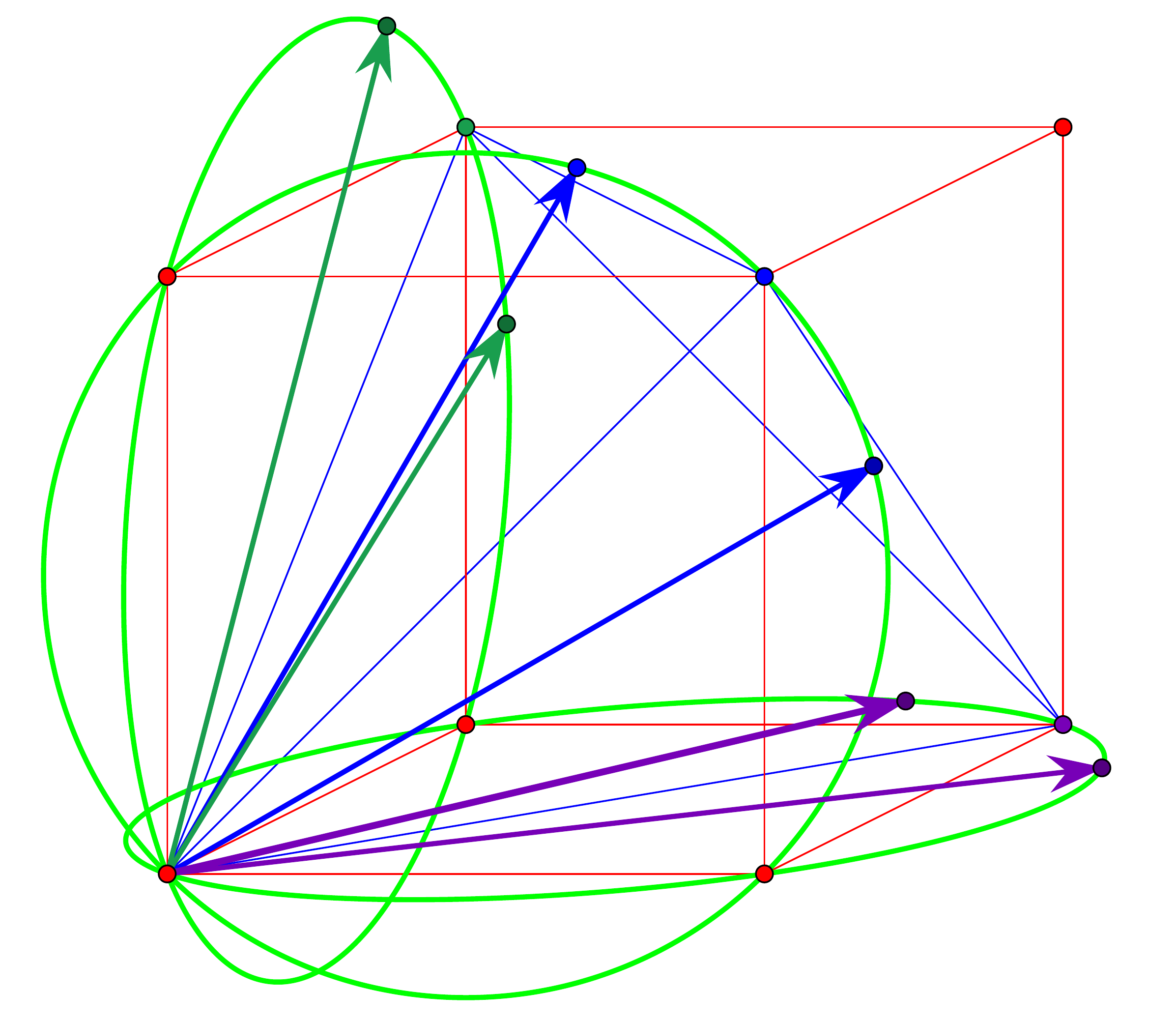}}
 \caption{The essentials of the diagram used for a three-dimensional framework with $n=4$ vertex orbits and $m=12$ edge orbits. The reduced (quotient) graph is represented by the regular tetrahedron inscribed in the cube. The six edge vectors emanating from vertex are shown as arrows. Only three great circles of the six spheres are traced. Tetrahedral symmetry completes the information recorded by the diagram.}
 \label{Fign4m12}
\end{figure}

The twenty-four transforms of this edge vector give representatives for the twelve edge orbits in twelve pairs with opposite signs. The six edge vectors shown in Figure~\ref{Fign4m12} are:

$$ 12A=(\alpha,\beta,0), \ \ \ 12B=(\beta,\alpha,0), $$
$$ 13A=(0,\alpha,\beta), \ \ \ 13B=(0,\beta,\alpha), $$
$$ 14A=(\beta,0,\alpha), \ \ \ 14B=(\alpha,0,\beta). $$

\noindent
With $IJA=-JIA$ and $IJB=-JIB$, the remaining six edge vector representatives are:

$$ 23A=(-\alpha,0,\beta), \ \ \ 23B=(-\beta,0,\alpha), $$
$$ 24A=(0,-\beta,\alpha), \ \ \ 24B=(0,-\alpha,\beta), $$
$$ 34A=(\beta,-\alpha,0), \ \ \ 34B=(\alpha,-\beta,0). $$

\noindent
The periodicity lattice is spanned by the tetrahedral group transforms of vectors of the form $12A-12B=(\alpha-\beta,\beta-\alpha,0)$ and $12A+23A+31A=(0,\beta-\alpha,0)$, which illustrates a period  associated to a cycle in the reduced graph. It follows that, regardless of the rationality or irrationality of $\alpha \neq \beta$, the resulting periodicity lattice is
$(\alpha-\beta){Z}^3$.

\begin{figure}[h]
\centering
 {\includegraphics[width=0.40\textwidth]{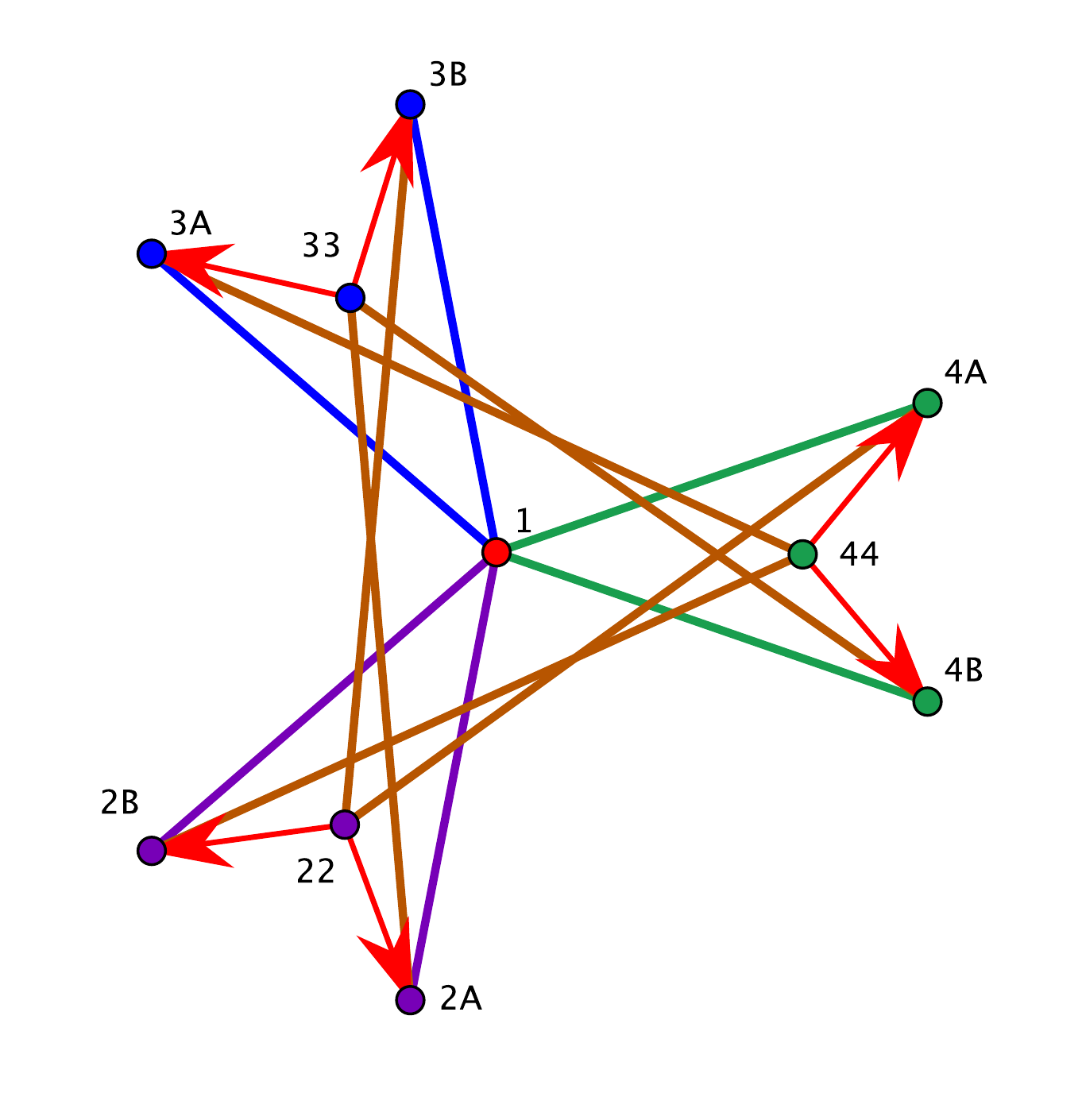}}
 \caption{The essentials of the resulting framework. The periodicity lattice
is generated by red vectors, with (22)(2A)=(44)(4B), (33)(3A)=(22)(2B), (44)(4A)=(33)(3B).}
 \label{FiginDepth}
\end{figure}

\medskip \noindent
The periodic framework itself is obtained by indefinite articulation of new edge
vectors, proceeding from the endpoints of the six depicted arrows and always respecting the vertex orbits. A partial extension, with all twelve edge representatives, is shown, from a different perspective, in Figure~\ref{FiginDepth}. With the vertex marked 1 placed at the origin (and seen ``in depth''), one may observe three skew quadrilaterals, 
resulting from the following relations.

$$ 12A+23A=14B+43B=(0,\beta,\beta)=33   $$
$$ 13A+34A=12B+24B=(\beta,0,\beta)=44   $$
$$ 14A+42A=13B+32B=(\beta,\beta,0)=22   $$

As a consequence of tetrahedral symmetry for the initial diagram, the resulting framework has a crystallographic symmetry group which
is transitive on vertices and transitive on edges. The one-parameter auxetic 
deformation which maintains this crystallographic symmetry and extends the strictly auxetic infinitesimal deformation
provided by our  construction can be described as follows. Dilate the (cube and inscribed) tetrahedron of the initial diagram, but maintain the squared length
of the new edge vector representative $(\tilde{\alpha}, \tilde{\beta},0)$, that is
$ \tilde{\alpha}^2 + \tilde{\beta}^2=\alpha^2+\beta^2$. Thus, at the vertex marked 1, the edge pairs $(12A,12B)$, $(13A,13B)$,
$(14A,14B)$ increase their angle in their respective planes (and by symmetry, this
type of local motion is replicated at any other vertex).

\medskip \noindent
{\bf Remark.}\ This family of periodic frameworks has (geometrically allowed)
self-intersections. It can be shown that these equal edge-length frameworks have four degrees of freedom i.e. a smooth four-dimensional local deformation space.
 
\section{Perusal of the virtual auxetic design catalog}

We have described above a procedure which associates to an \textit{initial diagram}
in $ \mathbb{R}^d$, a $d$-periodic framework with a strictly auxetic infinitesimal deformation. The initial diagram consists of the following elements: 

\medskip
(i) a finite connected multi-graph without loops, 

(ii) a placement of the (simple) reduced graph in $ \mathbb{R}^d$, 

(iii) spheres with diameters given by all edges of the placed reduced graph,

(iv) in a sphere, say, over the edge from $v_i$ to $v_j$, a depiction of (framework edge)
vectors $e_{ij}^k$, $k=1,...,m_{ij}$, where $m_{ij}$ is the number of edges in the 
multigraph over the specified edge in the reduced graph; all vectors $e_{ij}^k$
emanate from $v_i$ and have their endpoints on the sphere,

(v) a rank $d$ lattice generation condition for the span of the periodicity vectors; a periodicity vector is either a vector between the endpoints of
two (framework edge) vectors in the same sphere or a (cycle) vector of the form
$e_{i_1 i_2}^*+e_{i_2 i_3}+...+e_{i_J i_1}^*$ for a cycle of edges in the reduced graph.

\medskip
The associated periodic framework may turn out to be a singular point in its
deformation space and in this case additional investigations would be required
for deciding if the strictly auxetic infinitesimal deformation belongs to a
local auxetic path or not. If the periodic framework has a smooth local deformation space, the are always strictly auxetic path extensions.

\medskip
Thus, for a systematic and unencumbered generation of auxetic designs, we need a direct guarantee of property (v) for initial diagrams and a simple smoothness
criterion for the local deformation space of the resulting periodic framework.
We have already shown that, when operating with \textit{rational initial diagrams},
that is, diagrams which have the vertices of the reduced graph at rational points
and the endpoints of the (edge framework) vectors at rational points of the
corresponding spheres, property (v) is guaranteed to hold. For the smoothness issue, a sensible restriction is to use as multi-graphs only \textit{quotient graphs of minimally rigid periodic frameworks, with one or more edges removed}.
The class of quotient graphs of minimally rigid periodic frameworks has been
characterized in combinatorial terms in \cite{borcea:streinu:MinimallyRigid:BLMS:2011}.
With this restriction, smoothness is guaranteed by a single (maximal rank)
linear test.

\medskip \noindent
{\bf Classification outlook for two vertex orbits.}\ It is natural to organize the
virtual catalog according to increasing values of the number $n$ of vertex orbits
and the number $f$ of degrees of freedom. For $n=2$, the reduced (quotient) graph
is a single edge. Thus, we have a single ellipsoid and no cycles. We consider here
the case of periodic frameworks with $f=1$, i.e. one degree of freedom, for dimensions $d=2$ and $d=3$. 
For independent edge constraints, we have $f=dn+{d\choose 2}-m$, with
$m$ denoting the number of edge orbits.

\medskip
In the planar case, we have $m=4$. A periodic framework (with four independent
edge orbit constraints) is immediately recognized as a {\em strictly auxetic mechanism}
by tracing the unique conic through a vertex and the four endpoints of the bars
emanating from it: this conic must be an {\em ellipse}. Conversely, planar auxetic
periodic mechanisms can be generated by adopting as periodicity lattice $ \mathbb{Z}^2\subset \mathbb{R}^2$. We choose four points in $  \mathbb{Z}^2$, forming
a (non-degenerate) convex quadrilateral. In general, the $\mathbb{Z}$-span of the
edge vectors will be a sublattice of $  \mathbb{Z}^2$, but we shall retain only
choices with equality (without any loss of design patterns). For classification purposes, one must identify equivalent choices of this nature.
There is a one parameter family of ellipses running through four such points. We may select any ellipse and place the vertex representative of the second vertex orbit on it (away from the four points). 
The four (framework) edge vectors from this new point to the previous 
four points provide the data for articulating the associated auxetic periodic mechanism.

\medskip
In dimension three, we have $m=8$.  A periodic framework (with eight independent edge orbit constraints) is immediately recognized as a {\em strictly auxetic mechanism} by tracing the unique quadric surface through a vertex and the eight endpoints of the bars emanating from it: this quadric must be an {\em ellipsoid}. The design procedure is slightly more complicated than in the planar case because eight points in strictly convex position need not lie on any ellipsoid. We start with eight points in $\mathbb{Z}^3$ in strictly convex position and check that the $\mathbb{Z}$-span of the vectors between all point pairs equals $\mathbb{Z}^3$. Then we check that the family of quadrics passing through the eight points is one dimensional (i.e. a pencil) and contains ellipsoids. If this is the case, the choice is valid. Again, for classification purposes, one must identify equivalent choices of this nature. Auxetic mechanisms will be obtained, as above, by selecting an ellipsoid through the eight points and a vertex representative of the second vertex orbit on it (away from the quartic curve through the eight points which is the axis of the pencil). Up to an affine transformation, Figure~\ref{Cover} illustrates one such  blueprint.

\section{Conclusion}

Although auxetic behavior is typically non-linear, we have found a structural characterization of periodic auxetics from infinitesimal i.e. linear considerations. This necessary and sufficient condition for periodic bar-and-joint frameworks with strictly auxetic infinitesimal deformations is expressed in a diagram with homothetic ellipsoids over the reduced quotient graph and leads to endless, yet systematic possibilities for generating auxetic periodic designs.

\medskip \noindent
{\bf Acknowledgement.}\ This work was supported by the National Science Foundation (award no. 1319389 to C.S.B., award no. 1319366 to I.S.) and the National Institutes of Health (award no. 1R01GM109456 to C.S.B. and I.S.).


\begin{thebibliography}{10}

\bibitem{bertoldi:cellularReview:2017}
K.~Bertoldi.
\newblock Harnessing instabilities to design tunable architected cellular
  materials.
\newblock {\em Annu. Rev. Mater. Res.}, 47:51--61, 2017.

\bibitem{borcea:streinu:PeriodicFF:RSPA:2010}
Ciprian~S. Borcea and Ileana Streinu.
\newblock Periodic frameworks and flexibility.
\newblock {\em Proc R Soc Lond A Math Phys Sci}, 466:2633--2649, 2010.

\bibitem{borcea:streinu:MinimallyRigid:BLMS:2011}
Ciprian~S. Borcea and Ileana Streinu.
\newblock Minimally rigid periodic graphs.
\newblock {\em Bulletin of the London Mathematical Society}, 43:1093--1103,
  2011.

\bibitem{borcea:streinu:CrystallographicSym:PhT:2014}
Ciprian~S. Borcea and Ileana Streinu.
\newblock Frameworks with crystallographic symmetry.
\newblock {\em Philos Trans A Math Phys Eng Sci}, 372:20120143, 2014.

\bibitem{borcea:streinu:kinematicsExpansive:ark14:2014}
Ciprian~S. Borcea and Ileana Streinu.
\newblock Kinematics of expansive planar periodic mechanisms.
\newblock In Jadran Lenarcic and Oussama Khatib, editors, {\em Advances in
  Robot Kinematics (ARK'14)}, pages 395--408. Springer Verlag, Ljubljana,
  Slovenia, June 2014.

\bibitem{borcea:streinu:GeomAuxetics:RSPA:2015}
Ciprian~S. Borcea and Ileana Streinu.
\newblock Geometric auxetics.
\newblock {\em Proc R Soc Lond A Math Phys Sci}, 471:20150033, 2015.

\bibitem{borcea:streinu:LiftStress:DCG:2015}
Ciprian~S. Borcea and Ileana Streinu.
\newblock Liftings and stresses for planar periodic frameworks.
\newblock {\em Discrete Comput Geom}, 53:747--782, 2015.

\bibitem{borcea:streinu:AuxeticElliptic:2016}
Ciprian~S. Borcea and Ileana Streinu.
\newblock Auxetic deformations and elliptic curves.
\newblock {\em arXiv:1612.02100}, preprint, 2016.

\bibitem{borcea:streinu:Diamond:imaMath:2017}
Ciprian~S. Borcea and Ileana Streinu.
\newblock Deforming diamond.
\newblock {\em IMA J Appl Math}, 82(2):371--383, April 2017.

\bibitem{borcea:streinu:NewPrinciples:SIAGA:2017}
Ciprian~S. Borcea and Ileana Streinu.
\newblock New principles for auxetic periodic design.
\newblock {\em SIAM Journal on Applied Algebra and Geomatry}, 1(1):442--458, 2017.
\newblock DOI:10.1137/16M1088259.

\bibitem{dove:displacive:1997}
Martin~T. Dove.
\newblock Theory of displacive phase transitions in minerals.
\newblock {\em American Mineralogist}, 82:213--244, 1997.

\bibitem{elipe:lantada:auxeticGeometries:2012}
J.~C.~A. Elipe and A.~D. Lantada.
\newblock Comparative study of auxetic geometries by means of computer-aided
  design and engineering.
\newblock {\em Smart Materials and Structures}, 21:105004, 2012.

\bibitem{evans:auxetic:actaMetal:1994}
Kenneth~E. Evans, M.~A. Nkansah, and I.~J. Hutchinson.
\newblock Auxetic foams: Modeling negative {P}oisson's ratios.
\newblock {\em Acta Metall. Mater.}, 42(4):1289--1294, 1994.

\bibitem{evans:etAl:molecularNetwork:Nature:1991}
Kenneth~E. Evans, M.~A. Nkansah, I.~J. Hutchinson, and S.~C. Rogers.
\newblock Molecular network design.
\newblock {\em Nature}, 353:124--125, 1991.

\bibitem{greaves:lakes:etAl:PoissonRatio:2011}
G.~N. Greaves, A.~I. Greer, R.~Lakes, and T.~Rouxel.
\newblock Poisson's ratio and modern materials.
\newblock {\em Nature Materials}, 10:823--837, 2011.

\bibitem{grima:negativePoisson:RSPA:2012}
J.~N. Grima, R.~Caruana-Gauci, D.~Attard, and R.~Gatt.
\newblock Three-dimensional cellular structures with negative {P}oisson's ratio
  and negative compressibility properties.
\newblock {\em Proc R Soc Lond A Math Phys Sci}, 468:3121--3138, 2012.

\bibitem{huang:chen:negativePoisson:2016}
C.~Huang and L.~Chen.
\newblock Negative {P}oisson's ratio in modern functional materials.
\newblock {\em Advanced Materials}, 28:8079--8096, 2016.

\bibitem{kolken:auxetic:rscAdvances:2017}
H.~M.~A. Kolken and A.~A. Zadpoor.
\newblock Auxetic mechanical metamaterials.
\newblock {\em RSC Advances}, 7:5111--5129, 2017.

\bibitem{lakes:negativePoisson:1987}
R.~Lakes.
\newblock Foam structures with a negative {P}oisson's ratio.
\newblock {\em Science}, 235:1038--1040, 1987.

\bibitem{lakes:auxeticReview:2017}
R.~S. Lakes.
\newblock Negative-{P}oisson's-ratio materials: Auxetic solids.
\newblock {\em Annu. Rev. Mater. Res.}, 47:63--81, 2017.

\bibitem{lee:singer:thomas:microNanoMaterials:advMat:2012}
Jae-Hwang Lee, Jonathan~P. Singer, and Edwin~L. Thomas.
\newblock Micro-/nanostructured mechanical metamaterials.
\newblock {\em Advanced Materials}, 24:4782--4810, 2012.

\bibitem{mitschkeEtAl:geometryLeading:2016}
H.~Mitschke, F.~Schury, K.~Mecke, F.~Wein, M.~Stingl, and G.~E. Schroeder-Turk.
\newblock Geometry: The leading parameter for the {P}oisson's ratio of bending-
  dominated cellular solids.
\newblock {\em International Journal of Solids and Structures}, 10-101:1--10,
  2016.

\bibitem{park:auxetic:appliedPhysics:2016}
H.~S. Park, S.~Y. Kim, and J-W. Jiang.
\newblock Auxetic nanomaterials: Recent progress and future development.
\newblock {\em Applied Physics Reviews}, 3:041101, 2016.

\bibitem{reisEtAl:designerMatter:2015}
P.~M. Reis, H.~M. Jaeger, and M.~van Hecke.
\newblock Designer matter: A perspective.
\newblock {\em Extreme Mechanics Letters}, 5:25--29, 2015.

\bibitem{saxena:auxeticsResearch:advancedEngr:2016}
K.~K. Saxena, R.~Das, and E.~P. Calius.
\newblock Three decades of auxetics research: Materials with negative
  {P}oisson's ratio: A review.
\newblock {\em Advanced Engineering Materials}, 18(11):1847--1870, 2016.

\bibitem{lubensky:twistedKagome:2012}
K.~Sun, A.~Souslov, X.~Mao, and T.~C. Lubensky.
\newblock Surface phonons, elastic response, and conformal invariance in
  twisted kagome lattices.
\newblock {\em Proc. Nat. Acad. Sci.}, 109:12369--12374, 2012.

\end{thebibliography}
\end{document}